\documentclass[12pt]{article}
\usepackage{amssymb}
\title{Quadratic embeddings}
\author{Hans Havlicek \and Corrado Zanella}

\newtheorem{defi}{Definition}
\newtheorem{prop}{Proposition}[section]
\newtheorem{theo}{Theorem}
\newtheorem{lemma}{Lemma}
\newtheorem{example}{Example}

\newcommand{\Acal}{{\cal A}}
\newcommand{\Bcal}{{\cal B}}
\newcommand{\Ccal}{{\cal C}}
\newcommand{\Ecal}{{\cal E}}
\newcommand{\Fcal}{{\cal F}}
\newcommand{\Hcal}{{\cal H}}
\newcommand{\Kcal}{{\cal K}}
\newcommand{\Lcal}{{\cal L}}
\newcommand{\Mcal}{{\cal M}}

\newcommand{\Pcal}{{\cal P}}

\newcommand{\Scal}{{\cal S}}
\newcommand{\Tcal}{{\cal T}}
\newcommand{\Ucal}{{\cal U}}
\newcommand{\Xcal}{{\cal X}}

\newcommand{\CC} 
{\mbox{\sf \hspace{.25em}\rule[.05ex]{.1em}{1.5ex}%
\hspace{-.35em}C}}
\newcommand{\NN}{\mbox{\sf I\hspace*{-.12em}N}}

\newcommand{\ee}{\mbox{\bf e}}

\newcommand{\PG}[2]{\mbox{$\mbox{{\rm PG}}(#1,#2)$}}
\newcommand{\AG}[2]{\mbox{$\mbox{{\rm AG}}(#1,#2)$}}
\newcommand{\clos}[1]{\overline{\overline{#1}}}
\newcommand{\applicaz}[3]{\mbox{$#1\,:\,#2\rightarrow#3$}}
\newcommand{\Applicaz}[3]{\mbox{$#1\,:\,#2\longrightarrow#3$}}

\newcommand{\im}{\mbox{\rm im\,}}

\newcommand{\proof}{\noindent {\sl Proof}\hspace{.2in}}

\newcommand{\overnu}{{\widehat{\nu}}}

\newcommand{\citt}[4]{\bibitem{#1}
{\sc #2}: {\sl #3\/}. #4.}

\newcommand{\Ucek}{\check{\cal U}}

\begin{document}
\maketitle

   \begin{abstract}
   The quadratic Veronese embedding $\rho$ maps the point set $\Pcal$ of
   \PG{n}{F} into the point set of \PG{{n+2 \choose 2}-1}F ($F$ a commutative
   field) and has the following well-known property: If $\Mcal\subset\Pcal$,
   then the intersection of all quadrics containing $\Mcal$ is the inverse
   image of the linear closure of $\Mcal^{\rho}$. In other words, $\rho$
   transforms the closure from quadratic into linear. In this paper we use
   this property to define ``quadratic embeddings''. We shall prove that if
   $\nu$ is a quadratic embedding of $\PG{n}{F}$ into $\PG{n'}{F'}$ ($F$ a
   commutative field), then $\rho^{-1}\nu$ is dimension-preserving.
   Moreover, up to some exceptional cases, there is an injective homomorphism
   of $F$ into $F'$. An additional regularity property for quadratic
   embeddings allows us to give a geometric characterization of the
   quadratic Veronese embedding.
   \end{abstract}


\section{Introduction}

The aim of this paper is to examine geometric properties of a {\em
quadratic embedding}, i.e. a mapping between projective spaces sharing some
properties of the classical quadratic Veronese embedding. We follow an
approach
that has been used in discussing embeddings of Grassmann spaces (cf.\
\cite{Ha81}
and \cite{We83}) and product spaces  (cf.\ \cite{Za94}). See also
\cite[chapter 25]{HT91} for combinatorial characterizations of Veronese
varieties over finite fields.

Let $F$ be a commutative field and $(\Pcal,\Lcal):=\PG{n}{F}$. Write
   \begin{displaymath}
   \Phi:=\{\Scal\subset\Pcal|\,\Scal \mbox{ is a quadric of }
   \PG{n}{F}\}\cup\{\Pcal\}.
   \end{displaymath}
If $\Mcal\subset\Pcal$, then the {\em quadratic closure\/} of
$\Mcal$ is
   \begin{displaymath}
   \clos{\Mcal}:=\bigcap_{\Mcal\subset\Scal,\,\Scal\in\Phi}\Scal.
   \end{displaymath}We call $\Mcal$ a {\em closed set\/} if
   $\Mcal=\clos{\Mcal}$.
The linear closure of a set $\Mcal$ of points will be denoted by
$\overline{\Mcal}$. Each hyperplane of \PG{n}{F} is a quadric, namely a
repeated hyperplane. Hence $\Mcal\subset\clos\Mcal\subset\overline\Mcal$.

   \begin{defi}
   Let $(\Pcal,\Lcal):=\PG{n}{F}$ and $(\Pcal',\Lcal'):=\PG{n'}{F'}$, where
   the field $F$ is commutative. A mapping \applicaz{\nu}{\Pcal}{\Pcal'}
   is a {\em quadratic embedding} if
      \begin{equation}\label{D/1/1}
      \clos{\Mcal}=
      (\overline{\Mcal^{\nu}}){}^{\nu^{-1}}
      \mbox{ for all } \Mcal\subset\Pcal,
      \end{equation}
   and
      \begin{equation}
      \overline{\im\nu}=\Pcal'.
      \end{equation}
   \end{defi}

We give some examples of quadratic embeddings:
   \begin{example}{\rm
   The classical {\em quadratic Veronese embedding\/} $\rho$ is defined in
   the case $F'=F$, $n'=
   {n+2 \choose 2}-1$, by
      \begin{displaymath}
      F(x_0,\ldots,x_n) \stackrel{\rho}{\longmapsto}
      F(y_{ij})_{0\leq i\leq j\leq n},
      \mbox{ with } y_{ij}:= x_ix_j.
      \end{displaymath}
   There are many equivalent definitions. Cf., e.g., \cite{Bu61},
   \cite{He82}, \cite{Ka87}.
   }\end{example}
   \begin{example}{\rm
   Let $n'={n+2 \choose 2}-1$. If \applicaz{\alpha}F{F'}\ is an injective
   homomorphism, then $\alpha$ induces a canonical embedding $\epsilon$ of
   \PG{n'}F\ into \PG{n'}{F'}. The mapping $\rho\epsilon$ turns out to be a
   quadratic embedding. Since there are examples of fields admitting an
   injective, but not surjective homomorphism \applicaz{\alpha}FF,
   there exist quadratic embeddings different from the classical one, even
   if we demand that $F$ and $F'$ are isomorphic.
   }\end{example}
   \begin{example}{\rm
   If $n=1$ then $\Mcal\subset\Pcal$ is closed if, and only if, $|\Mcal|\leq
   2$ or $\Mcal=\Pcal$. Thus a mapping \applicaz{\nu}{\Pcal}{\Pcal'} is a
   quadratic embedding if, and only if, $n'=2$, $\nu$ is injective and
   $\im\nu$ is an arc.
   }\end{example}
   \begin{example}{\rm
   If \PG nF = \PG 22 and $\Fcal$ is a frame in \PG 5{F'} then any injection
   \applicaz{\nu}{\Pcal}{\Pcal'} such that $\im\nu=\Fcal$ is a quadratic
   embedding. This is immediate from the fact that each subset $\Mcal$ of
   $\Pcal$ is closed, unless $|\Mcal|=6$.
   }\end{example}


\section{Properties of quadratic embeddings}

In this section $\nu$ is a quadratic embedding of $\PG nF=(\Pcal,\Lcal)$
($n\geq 1$) into $\PG{n'}{F'}=(\Pcal',\Lcal')$.

   \begin{prop}\label{3/5}
   Let $\Kcal_1$ and $\Kcal_2$ be two distinct closed sets in \PG nF. Then
   $\overline{\Kcal^{\nu}_1}\neq \overline{\Kcal^{\nu}_2}$. Consequently,
   the mapping $\nu$ is injective and satisfies
      \begin{equation}\label{3/5/1}
      \clos{\Mcal}^{\nu}=\overline{\Mcal^{\nu}}\cap \im\nu
      \mbox{ for all }
      \Mcal\subset\Pcal.
      \end{equation}
   If\/ $\Ucal'\subset\Pcal'$ is a subspace, then
   $\Ucal'^{\nu^{-1}}\subset\Pcal$ is a closed set.
   \end{prop}
\proof By the definition of a quadratic embedding,
   \begin{displaymath}
   (\overline{\Kcal^{\nu}_1})^{\nu^{-1}}=\clos{\Kcal_1}=\Kcal_1\not=
   \Kcal_2=\clos{\Kcal_2}=(\overline{\Kcal^{\nu}_2})^{\nu^{-1}},
   \end{displaymath}
so that $\overline{\Kcal^{\nu}_1}\neq \overline{\Kcal^{\nu}_2}$. Moreover,
$\nu$ is injective, since any subset of $\Pcal$ with a single element is
closed. Hence (\ref{3/5/1}) is true. Finally, let $\Mcal:=\Ucal'^{\nu^{-1}}$.
Then
   \begin{displaymath}
   \clos{\Mcal}=(\overline{\Mcal^{\nu}}){}^{\nu^{-1}}
   \subset\Ucal'^{\nu^{-1}}=\Mcal.\Box
   \end{displaymath}

   \begin{theo}\label{3/7}
   If $\nu$ is a quadratic embedding of $\PG nF$
   into $\PG{n'}{F'}$,
   then
   $n'={n+2 \choose 2}
   -1.$
   \end{theo}

\proof Define $\delta(t):=
{t+2 \choose 2}$, $t\in\NN$. Let $\{\ee_0,\ldots,\ee_n\}$ be a basis of
$F^{n+1}$ and
   \begin{displaymath}\Xcal:=\{F(\ee_i+\ee_j)|0\leq i<j\leq n\}
   \cup\{F\ee_i|i=0,\ldots,n\}.
   \end{displaymath}
Since $|\Xcal|=\delta(n)$ and $\clos{\Xcal}=\Pcal$, we have
$\im\nu\subset\overline{\Xcal^{\nu}}$, hence
   \begin{equation}\label{3/7/1}
   n'\leq\delta(n)-1.
   \end{equation}

We now prove that in (\ref{3/7/1}) the equality holds. We give a definition,
by recursion on $d=0,\ldots,n$, of distinct closed sets in \PG nF, say
$\Kcal_{\delta(d-1)}$, $\Kcal_{\delta(d-1)+1}$, \ldots,
$\Kcal_{\delta(d)-1}$, such that
   \begin{displaymath}
   \Kcal_{\delta(d-1)}\subset\Kcal_{\delta(d-1)+1}\subset
   \ldots\subset\Kcal_{\delta(d)-1},
   \end{displaymath}
with $^d\Ucal:=\Kcal_{\delta(d)-1}$ being a $d$-subspace of \PG nF. For
$d=0$, choose a point $Q$ and let $\Kcal_0={}^0\Ucal:=\{Q\}$. Now let $d>0$
and
$\Kcal_{\delta(d-1)-1}={}^{d-1}\Ucal$. Take a $d$-subspace $^d\Ucal$
containing $^{d-1}\Ucal$ and a basis $\Bcal=\{P_0,\ldots,P_d\}$ of $^d\Ucal$
such that $^{d-1}\Ucal\cap\Bcal=\emptyset$. Since the union of two subspaces
of \PG nF\ is a closed set, we can define
   \begin{displaymath}
   \Kcal_{\delta(d-1)+i}:={}^{d-1}\Ucal\cup \overline{\{P_0,\ldots,P_i\}},
   \hspace{.5in}i\in\{0,\ldots,d\}.
   \end{displaymath}
By Prop.~\ref{3/5},
\begin{displaymath}
   \emptyset\subset \overline{\Kcal^{\nu}_0}\subset
   \overline{\Kcal^{\nu}_1}\subset\ldots\subset \Kcal_{\delta(n)-1}^{\nu}
   \end{displaymath}
is a chain of distinct subspaces of \PG{n'}{F'}.$\Box$

   \begin{prop}\label{3/9}
   If $\Mcal\subset\Pcal$, then $\dim(\overline{\Mcal^{\nu}})$ is equal to
   the largest $i\in\NN$, such that there exists a chain
      \begin{equation}\label{3/9/1}
   \emptyset\subset\Kcal_0\subset\Kcal_1\subset\ldots
      \subset\Kcal_i=\clos{\Mcal},
      \end{equation}
   consisting of $i+2$ distinct closed subsets of $\clos{\Mcal}$.
   Consequently,
   $\dim(\overline{\Mcal^{\nu}})=\dim(\overline{\Mcal^{\rho}})$,
   where $\rho$ denotes the quadratic Veronese embedding.
   \end{prop}

\proof By Prop.~\ref{3/5}, the subspaces $\overline{\Kcal^{\nu}_0}$,
$\overline{\Kcal^{\nu}_1}$, \ldots, $\overline{\Kcal^{\nu}_i}$ are distinct
and
   \begin{displaymath}
   \overline{\Kcal^{\nu}_i}= \overline{\overline{\Mcal^{\nu}}\cap\im\nu}=
   \overline{\Mcal^{\nu}},
   \end{displaymath}
whence $\dim(\overline{\Mcal^{\nu}})\geq i$.

Now assume $\dim(\overline{\Mcal^{\nu}})>i$. Then there exists an integer
$j$, $0\leq j<i$, such that
   \begin{displaymath}
   \dim(\overline{\Kcal^{\nu}_{j+1}})\neq \dim(\overline{\Kcal^{\nu}_j})+1.
   \end{displaymath}
Let $P\in\Kcal^{\nu}_{j+1}\setminus\overline{\Kcal^{\nu} _j}$. Then
$\Kcal:=(\overline{\{P\}\cup\Kcal^{\nu}_j}){}^{\nu^{-1}}$ is a closed set
(cf.\ Prop.~\ref{3/5}), and $\Kcal_j\subset\Kcal\subset\Kcal_{j+1}$,
$\Kcal\neq\Kcal_j$. The maximality of the chain (\ref{3/9/1}) implies
$\Kcal=\Kcal_{j+1}$. Therefore
   \begin{displaymath}
   \dim(\overline{\Kcal^{\nu}_{j+1}})=
   \dim(\overline{\{P\}\cup\Kcal^{\nu}_j})= \dim(\overline{\Kcal^{\nu}_j})+1,
   \end{displaymath}
a contradiction.$\Box$

   \begin{prop}\label{immsing}
   If $\Tcal$ is a hyperplane of \PG nF\ and $\Mcal\subset \Pcal\setminus
   \Tcal$, then
      \begin{equation}\label{immsing/1}
      \dim(\overline{(\Tcal\cup\Mcal)^{\nu}})=
         {n+1 \choose 2}
         +\dim(\overline{\Mcal}).
      \end{equation}
   \end{prop}

\proof By Theorem~\ref{3/7}, $\dim(\overline {\Tcal^{\nu}})=
         {n+1 \choose 2}-1$. If
$\Bcal=\{P_0,\ldots,P_t\}\subset\Mcal$ is a basis of $\overline{\Mcal}$, then
the closed sets
   \begin{displaymath}
   \Kcal_i:=\Tcal\cup\overline{\{P_0,\ldots,P_i\}},
   \hspace{.5 in}i\in\{0,\ldots,t\},
   \end{displaymath}
form a saturated chain
   \begin{displaymath}
   \Tcal\subset\Kcal_0\subset\ldots\subset\Kcal_t.
   \end{displaymath}
Now the assertion is  a consequence of Theorem~\ref{3/7}
and Prop.~\ref{3/9}.$\Box$

   \begin{prop}\label{H/2}
   Let $\Tcal$ be a hyperplane of \PG nF. If
   $\overline{\Tcal^{\nu}}$ and $\Ecal'$ are complementary subspaces of
   \PG{n'}{F'}, then the mapping
      \begin{equation}\label{H/2/1}
      \Applicaz{\iota}{\Pcal\setminus \Tcal}{\Ecal'}
      \,:\, A\longmapsto \overline{(\Tcal\cup\{A\})^{\nu}}\cap\Ecal'
      \end{equation}
   has the following property:
      \begin{equation}\label{H/2/2}
      \dim(\overline{\Mcal})=\dim(\overline{\Mcal^{\iota}})
      \mbox{ for all }
      \Mcal\subset\Pcal\setminus \Tcal.
      \end{equation}
   Consequently, $\iota$ is preserving both collinearity and
   non-collinearity of points. So, the mapping $\iota$ is a {\em (linear)
   embedding} of the affine space $\Pcal\setminus \Tcal$ into the projective
   space $\Ecal'$.
   \end{prop}

\proof By (\ref{immsing/1}), $\dim\overline{\Tcal^{\nu}}=
{n+1 \choose 2}-1$, whence
$\dim\Ecal'=n$. Applying (\ref{immsing/1}) again yields
   \begin{eqnarray*}
   \dim(\overline{\Mcal^{\iota}}) &=&
   \dim(\overline{\Tcal^{\nu}\cup\Mcal^{\nu}}\cap\Ecal')=
   \dim(\overline{\Tcal^{\nu}\cup\Mcal^{\nu}})-
   (\dim(\overline{\Tcal^{\nu}})+1) \\
   &=& \dim(\overline{\Mcal}).\Box
   \end{eqnarray*}

   \begin{prop}\label{H/3}
   Let $|F|>2$ and $n\geq2$. If $|F|\not=3$ or $n\not=2$, then the embedding
   (\ref{H/2/1}) can be extended to exactly one embedding
   \applicaz{\beta}{\Pcal}{\Ecal'}.
   \end{prop}
\proof The case $n=2$ is dealt with in \cite{Ri65}. The case $|F|>3$ is
covered by \cite[Theorem 3.5]{BR84}. Thus only $|F|=3$ and $n>2$ remains
open. For $F'$ being finite,
the assertion follows from a result in \cite[chapitre
2.3]{Li80} (cf.\ also \cite[th\'eor\`eme 1]{Li82}), and by slight
modifications, this carries over to an infinite $F'$.

On the other hand we sketch a direct proof for $n>2$: Let $P\in \Tcal$. If
$g,h\in\Lcal$ and $P\in g\cap h$, then the lines
$\overline{(g\setminus\{P\})^{\iota}}$ and
$\overline{(h\setminus\{P\})^{\iota}}$ are coplanar by Prop.~\ref{H/2}.
Since
there exist three non coplanar lines through $P$, all lines of the kind
$\overline{(g\setminus\{P\})^{\iota}}$, with $P\in g$, share one point
$P'\in\Ecal'$. Then we define $P^{\beta}:=P'$. By repeatedly using
Prop.~\ref{H/2}, we have that $\beta$ is an embedding. The restriction
of $\beta$ to a plane $\Acal$ of $\Pcal$, not contained in $\Tcal$, is an
extension of $\iota|(\Acal\setminus\Tcal)$, and thus we obtain the
uniqueness of $\beta$.$\Box$

As a corollary, we have:
   \begin{theo}
   Let $|F|>2$ and $n\geq2$. If $|F|\not=3$ or $n\not=2$, then the existence
of
   a quadratic embedding of $\PG nF$ into $\PG{n'}{F'}$ implies that the
   field $F$ is isomorphic to a subfield of $F'$.$\Box$
   \end{theo}

Whenever for some fixed hyperplane $\Tcal\subset\Pcal$ and an adequately
chosen subspace $\Ecal'\subset\Pcal'$ the mapping (\ref{H/2/1}) is uniquely
extendable to an embedding \applicaz{\beta}{\Pcal}{\Ecal'}, then $\Tcal$
gives rise to an embedding
   \begin{equation}\label{nu_t}
   \Applicaz{\nu_\Tcal}{\Pcal}{\Pcal'/\overline{\Tcal^\nu}}\,:\,
   X\longmapsto\{X^\beta\}\vee\overline{\Tcal^\nu};
   \end{equation}
here $\Pcal'/\overline{\Tcal^\nu}$ denotes the point set of the quotient
space \PG{n'}{F'} modulo $\overline{\Tcal^\nu}$. Moreover, we can associate
with $\Tcal$ the following hyperplane of \PG{n'}{F'}:
   \begin{displaymath}
   \overline{\Tcal^{\nu}\cup \Tcal^{\beta}}=:\Hcal'_\Tcal.
   \end{displaymath}
Both definitions do not depend on the choice of $\Ecal'$. Since
$\Hcal'_\Tcal\cap\Ecal'\cap\im\iota=\emptyset$, we have
   \begin{equation}\label{**}
   (\Hcal'_\Tcal)^{\nu^{-1}}=\Tcal.
   \end{equation}

   \begin{prop}\label{H/4}
   Let $\Sigma$\ be the collection of all hyperplanes of \PG nF. If
   $\Hcal'_\Tcal$ is well defined for all $\Tcal\in\Sigma$, then
   \begin{equation} \label{overnu}
        \Applicaz{\overnu}{\Sigma}{\Sigma'}\,:\, \Tcal\longmapsto\Hcal'_\Tcal
   \end{equation}
   is an injective mapping.$\Box$
   \end{prop}

The previous results give sufficient conditions for the existence of
the mapping $\overnu$.


\section {Regular quadratic embeddings}

In the following we shall assume that $\nu$ is a quadratic embedding of
\PG nF, $n\geq 1$, into \PG{n'}{F'}, $n'={n+2 \choose 2}-1$.

\begin{defi}
  A quadratic embedding $\nu$ is called {\em $(P,\ell)$-regular\/} if there
  exists an incident point-line pair $(P,\ell)$ of \PG nF\ such that the
  plane arc $\ell^{\nu}$ has a unique unisecant line which is running
  through $P^{\nu}$ and contained in the plane $\overline{\ell^{\nu}}$.
  If $\nu$ is $(P,\ell)$-regular for all incident pairs $(P,\ell)$, then
  $\nu$ is said to be a {\em regular\/} quadratic embedding.
\end{defi}

\begin{prop} \label{X31}
  Suppose that $n\geq2$ and that $\nu$ is $(P,\ell)$-regular.
  Then $\nu$ is regular.
\end{prop}

\proof
By $n\geq2$, there exists a hyperplane $\Tcal$ of \PG nF\ such that
$P\in\Tcal$, $\ell\not\subset\Tcal$. Define an embedding
\applicaz{\iota}{\Pcal\setminus\Tcal}{\Ecal'} according to
(\ref{H/2/1}).
The $(P,\ell)$-regularity of $\nu$ implies that
\begin{equation} \label{XA}
  |\overline{(\ell\setminus\{P\})^{\iota}}\setminus
  (\ell\setminus\{P\})^{\iota}|=1.
\end{equation}

If the settings of Prop.~\ref{H/3} are true, then $\iota$ extends to an
embedding \applicaz{\beta}{\Pcal}{\Ecal'} with
$\ell^{\beta}=\overline{(\ell\setminus\{P\})^{\iota}}$ by (\ref{XA}).
Hence $\beta$ is a collineation.

Otherwise $|F|=:p\in\{2,3\}$ so that $|F'|=p$ by (\ref{XA}).
If $X\in\Pcal\setminus\Tcal$, then there is  a certain number of lines
through $X$ and on each such line there are $p$ points of
$\Pcal\setminus\Tcal$.
In $\Ecal'$ the same number of lines is running through $X^{\iota}$ and,
by Prop.~\ref{H/2}, there are $p$ points of $\im\iota$ on each such line.
This in turn means that on each line in $\Ecal'$ through a point of
$\Ecal'\setminus\im\iota$ there is either no point of $\im\iota$
or no other point of $\Ecal'\setminus\im\iota$, whence
$\Ecal'\setminus\im\iota$ is a subspace.
More precisely, $\Tcal':=\Ecal'\setminus\im\iota$ is a hyperplane of
$\Ecal'$.
Two distinct lines of the affine space $\Pcal\setminus\Tcal$ are parallel if,
and only if, they are disjoint and coplanar.
By Prop.~\ref{H/2} these properties carry over to the $\iota$-images of these
lines, whence $\iota$ is an affinity of $\Pcal\setminus\Tcal$ onto
$\Ecal'\setminus\Tcal'$.
(This is trivial when $p=3$.)
Thus $\iota$ is also extendable to a collineation
\applicaz{\beta}{\Pcal}{\Ecal'} if Prop.~\ref{H/3} cannot be applied.

Next choose any point $X_1\in\Tcal$ and any line
$\ell_1\not\subset\Tcal$, $X_1\in\ell_1$.
Then
$(\{X_1^\beta\}\vee\overline{\Tcal^\nu})\cap\overline{\ell^{\nu}_1}$
is the only unisecant of
$\ell_1^{\nu}$ at $X_1^{\nu}$ within the plane
$\overline{\ell_1^{\nu}}$, since
$|\overline{(\ell_1\setminus\{X_1\})^{\iota}}
\setminus(\ell_1\setminus\{X_1\})^{\iota}|=1$.
Hence $\nu$ is $(X_1,\ell_1)$-regular.
Repeatedly using this last idea yields that $\nu$ is regular.$\Box$

As an immediate consequence of the proof of Prop.~\ref{X31} we have
\begin{prop} \label{X32}
  Let $\nu$ be a regular quadratic embedding and $n\geq2$.
  Choose any hyperplane $\Tcal\subset\Pcal$.
  Then the embedding \applicaz{\iota}{\Pcal\setminus\Tcal}{\Ecal'},
  defined according to (\ref{H/2/1}), is extendable to a unique collineation
  \applicaz{\beta}{\Pcal}{\Ecal'}.
  Consequently, $F$ and $F'$ are isomorphic fields, and
  \applicaz{\nu_{\Tcal}}{\Pcal}{\Pcal'/\overline{\Tcal^{\nu}}}
  (cf.~(\ref{nu_t})) is a collineation.$\Box$
\end{prop}

If $n=1$, then the $(P,\ell)$-regularity of $\nu$ implies
\[
  |F|=|\ell\setminus\{P\}|=|\ell^{\nu}\setminus\{P^{\nu}\}|=|F'|.
\]
This does not imply, however, that $\nu$ is regular.
If $n=1$ and $\nu$ is regular, then $\im\nu$ obviously is an oval
but not necessarily a conic.
Cf.~\cite{Bu79,BHL80,St95} for topological conditions that force an oval
to be a conic.

These results can be improved if we assume that $|F|=:q$ is finite.
Then $n=1$ and $\nu$ being $(P,\ell)$-regular yield
$|F|=|F'|=q$ so that $\im\nu$ is a $(q+1)$-arc in
\PG 2{F'}${}\cong{}$\PG 2q.
Hence $\im\nu$ is an oval, which in turn shows that $\nu$ is regular.
Moreover, by Segre's theorem, $\im\nu$ is a (regular) conic if $q$ is
odd; the last result is also true when $q\in\{2,4\}$.

The case $n=1$ is excluded from our further discussions. Hence
we may assume without loss of generality that $F=F'$, by
virtue of Prop.~\ref{X32}.

The following result will be used in order to characterize the $\nu$-images
of lines.

   \begin{lemma}\label{H/6}
   Let $P_0$ and $P_2$ be two distinct points of \PG2F and let $\sigma$ be a
   collineation of \PG2F taking $P_0$ to $P_2$, but not fixing the line
   $\overline{P_0P_2}$. Then
      \begin{displaymath}
      \Ccal :=\{X|\{X\}=x\cap x^\sigma,x\mbox{ is a line through }P_0 \}
      \end{displaymath}
   is containing three distinct collinear points if, and only if, $\sigma$ is
   a non-projective collineation.
   \end{lemma}
\proof If $\sigma$ is projective, then $\Ccal$ is a regular conic, whence it
does not contain three distinct collinear points.

If $\sigma$ is not projective, then let $\alpha\in \mbox{Aut}(F)$ be the
companion automorphism of $\sigma$. Set
   \begin{displaymath}
   \{P_1\}:=(\overline{P_0P_2})^{\sigma^{-1}}\cap
   (\overline{P_0P_2})^{\sigma}.
   \end{displaymath}
Choose some line $e$ running through $P_0$ but not containing $P_1$ or $P_2$
and define $\{E\}:=e\cap e^\sigma$. Then $(P_0,P_1,P_2,E)$ is an ordered
quadrangle; we may assume that this is the standard frame of reference. A
straightforward calculation yields
   \begin{displaymath}
   \Ccal =\{F(u_0 u_0^\alpha,u_0 u_1^\alpha,u_1 u_1^\alpha)|
   (0,0)\not=(u_0,u_1)\in F^2\}.
   \end{displaymath}
By $\alpha\not=\mbox{id}_F$, there exists an element $c\in F$ with
$c\not=c^\alpha$. Define $v\in F$ via $c=v^{\alpha\alpha}c^\alpha$. Thus
$v\not=0,1$ and $F(1,1,1)$, $F(1,v^\alpha, vv^\alpha)$ are distinct points of
$\Ccal$. With
      \begin{displaymath}
      w:=\frac{1+vv^\alpha c}{1+v^\alpha c}
      \end{displaymath}
we obtain
      \begin{displaymath}
      \frac{1}{1+c}(1,1,1) + \frac{c}{1+c}(1,v^\alpha, vv^\alpha) =
      (1,w^\alpha,ww^\alpha),
      \end{displaymath}
whence $\Ccal$ is containing three distinct collinear points.$\Box$

\begin{prop}\label{H/7}
   If $g$ is a line of \PG nF, $n\geq2$, and $\nu$ is a regular quadratic
   embedding, then $g^\nu$ is a regular conic.
\end{prop}
\proof Choose
hyperplanes $\Tcal,\Ucek\subset\Pcal$ such that
$g\cap\Tcal\cap\Ucek=\emptyset$. Set
$\{T\}:=g\cap\Tcal$ and define a collineation $\nu_{\Tcal}$ according to
(\ref{nu_t}).
Write
   \begin{displaymath}
   \Lcal_T':=\{x'\in\Lcal'|T^{\nu}\in x'\subset\overline{g^\nu}\}
   \end{displaymath}
and
   \begin{displaymath}
   \Applicaz{\pi_T}{\Lcal_T'}{g^{\nu_{\Tcal}}}\,:\,
   x'\longmapsto x'\vee\overline{\Tcal^\nu}.
   \end{displaymath}
This $\pi_T$ is a projectivity from a pencil of lines onto a pencil of
subspaces.
Replacing $\Tcal$ by $\Ucek$ gives a point $U$ and a projectivity $\pi_U$.
Since $\nu_{\Tcal}^{-1}\nu_{\Ucek}$ is a collineation of quotient spaces,
   \begin{displaymath}
   \Applicaz{\pi_T\nu_{\Tcal}^{-1}\nu_{\Ucek}\pi_U^{-1}}{\Lcal_T'}{\Lcal_U'}
   \end{displaymath}
is extendable to a collineation, say $\sigma$, of the plane
$\overline{g^\nu}$ onto itself. We have
   \begin{displaymath}
      \begin{array}
      { r@{\:\stackrel{\sigma}{\longmapsto}\:} l }
      \overline{g^\nu}\cap \Tcal^\overnu   &\overline{T^\nu U^\nu}
      \not=\overline{g^\nu}\cap \Tcal^\overnu,
      \\
      \overline{T^\nu U^\nu}           &\overline{g^\nu}\cap \Ucek^\overnu,
      \\
      \overline{T^\nu X^\nu}           &\overline{U^\nu X^\nu}\,
      (X\in g\setminus\{T,U\})
      \end{array}
   \end{displaymath}
and
   \begin{displaymath}
   g^\nu=\{X'|\{X'\}=x'\cap x'^\sigma,x'\in\Lcal_T'\}.
   \end{displaymath}
Since any two distinct points of $g$ form a closed set, no three points of
$g^\nu$ are collinear. We read off from Lemma~\ref{H/6} that $\sigma$
is projective, whence $g^\nu$ is a regular conic.$\Box$

We remark that $\nu^{-1}_{\Tcal}\nu_{\Ucek}$ is a projective collineation
of quotient spaces.
\begin{prop} \label{X33}
  Let $(P_0,P_1,\ldots,P_n,E)$ be an ordered frame of \PG nF, $n\geq2$,
  and let $\nu$ be a regular quadratic embedding.
  Write $Q'_{ii}:=P^{\nu}_i$, $E':=E^{\nu}$, and $Q'_{ij}$ for the common
  point of the tangent lines of the conic $(\overline{P_iP_j})^{\nu}$
  at $P_i^{\nu}$ and $P_j^{\nu}$, $i,j\in\{0,1,\ldots,n\}$, $i\neq j$.
  Then
  \begin{equation}\label{frame}
     \{Q'_{ij}|0\leq i\leq j\leq n\}\cup\{E'\}
  \end{equation}
  is a frame of \PG {n'}F.
\end{prop}

\proof
For any $i,j\in\{0,1,\ldots,n\}$, $i<j$, take a
$P_{ij}\in\overline{P_iP_j}\setminus\{P_i,P_j\}$.
Then $\{Q'_{00},Q'_{11},\ldots,Q'_{nn}\}\cup\{P_{ij}^{\nu}|0\leq i<j\leq n\}$
is a basis of \PG {n'}F by Theorem \ref{3/7}.
Since $Q'_{ii}$, $Q'_{jj}$, $P_{ij}^{\nu}$, $Q'_{ij}$ is a plane quadrangle,
the exchange lemma yields that
\[
  \Bcal':=\{Q'_{ij}|0\leq i\leq j\leq n\}
\]
is a basis of \PG {n'}F.

Define hyperplanes
\[
  \Tcal_i:=\overline{\{P_k|k\in\{0,1,\ldots,n\}\setminus\{i\}\}}\subset\Pcal,
\]
and $\Xcal'_i:=\Tcal_i^{\overnu}$ (cf.~(\ref{overnu})) for
$i\in\{0,1,\ldots,n\}$.
Obviously $Q'_{jk}\in\Xcal'_i$ for all
$j,k\in\{0,1,\ldots,n\}\setminus\{i\}$.
Moreover, if $j\in\{0,1,\ldots,n\}\setminus\{i\}$, then
\[
  \overline{\{Q'_{ii},Q'_{ij},Q'_{jj}\}}\cap \Xcal'_i
\]
is the tangent line of the conic $(\overline{P_iP_j})^\nu$
at $P^{\nu}_j=Q'_{jj}$, so that $Q'_{ij}\in\Xcal'_i$.
We infer that
\[
  \Xcal'_i=\overline{\Bcal'\setminus\{Q'_{ii}\}}.
\]
Now $E\not\in\Tcal_i$ implies $E'\not\in\Xcal'_i$.
Finally, $\Tcal_i\cup\Tcal_j$ is a closed set not containing $E$.
Hence
\[
  E'\not\in\overline{(\Tcal_i\cup\Tcal_j)^\nu}=
  \overline{\Bcal'\setminus\{Q'_{ij}\}}.
\]
This completes the proof.$\Box$
\begin{theo}
If $\nu$ is a regular quadratic embedding of \PG nF into \PG{n'}F,
$n\geq2$, $n'={n+2\choose 2}-1$ and $\rho$ denotes the quadratic Veronese
embedding, then there exists a collineation $\kappa$ of \PG{n'}F such
that $\nu=\rho\kappa$.
\end{theo}

\proof
We adopt the notation of Prop.~\ref{X33}.
The coordinates with respect to $(P_0,P_1,\ldots,P_n,E)$ of a point
$X\in\Pcal$ are written as
$F(x_0,x_1,\ldots,x_n)$, and the coordinates of $X^\nu$ with respect
to $(Q'_{00},Q'_{01},\ldots,Q'_{nn},E')$ (cf.\ (\ref{frame}))
are denoted by
$F(y_{00},y_{01},\ldots,y_{nn})$.
In order to simplify notation we put $y_{ij}:=y_{ji}$ for $i>j$.

Choose an index $i\in\{0,1,\ldots,n\}$ and set
\[
  \Ecal'_i:=\overline{\{Q'_{i0},Q'_{i1},\ldots,Q'_{in}\}}.
\]
Hence $\Ecal'_i$ is a complement of $\overline{\Tcal^{\nu}_i}$
(cf.\ the proof of Prop.~\ref{X33}) and, by Prop.~\ref{X32}, we obtain a
collineation \applicaz{\beta_i}{\Pcal}{\Ecal'_i} with
$P_j\mapsto Q'_{ij}$ ($j\in\{0,1,\ldots,n\}$) and $E\mapsto E'_i$,
where $\{E'_i\}:=(\{E'\}\vee \overline{\Tcal_i^\nu})\cap\Ecal'_i$.
So, by taking $(Q'_{i0},Q'_{i1},\ldots,Q'_{in},E'_i)$ as frame of
reference in $\Ecal'$, we obtain that $X^{\beta_i}$ has coordinates
\[
  F(x_0^{\alpha_i},x_1^{\alpha_i},\ldots,x_n^{\alpha_i})
\]
with $\alpha_i\in\mbox{Aut}(F)$.

If $j\in\{0,1,\ldots,n\}$, then
$\nu^{-1}_{\Tcal_i}\nu_{\Tcal_j}$ is a projective collineation,
as has been remarked after Prop.~\ref{H/7}.
Hence also $\beta^{-1}_i\beta_j$ is projective.
Thus $\beta_i$ and $\beta_j$ belong to the same automorphism
$\alpha:=\alpha_i\in\mbox{Aut}(F)$.

Now we compare the coordinates of $X$, $X^\nu$, $X^{\beta_i}$:
If $x_i\neq0$, then
$\{X^{\beta_i}\}=(\{X^\nu\}\vee\overline{\Tcal_i^\nu})\cap\Ecal'_i$.
Hence there exists an element $c_i\in F\setminus\{0\}$ such that
\[
  y_{i0}=c_ix^{\alpha}_0,\
  y_{i1}=c_ix^{\alpha}_1,\ \ldots\ ,\
  y_{in}=c_ix^{\alpha}_n.
\]
If, moreover, $x_j\neq0$, $j\in\{0,1,\ldots,n\}\setminus\{i\}$, then
\[
  c_i=\frac{y_{ij}}{x^{\alpha}_j},\hspace{2em}
  c_j=\frac{y_{ji}}{x^{\alpha}_i}
\]
whence, by $y_{ij}=y_{ji}$,
\[
  \frac{c_i}{c_j}=\frac{x^{\alpha}_i}{x^{\alpha}_j}.
\]
If $x_i=0$, then $y_{i0}=y_{i1}=\ldots=y_{in}=0$.
Thus we have
\[
  F(y_{00},y_{01},\ldots,y_{nn})=
  F(x^{\alpha}_0x^{\alpha}_0,x^{\alpha}_0x^{\alpha}_1,\ldots,
  x^{\alpha}_nx^{\alpha}_n).
\]
Now letting $\kappa$ be that collineation of \PG{n'}F which transforms
each coordinate under $\alpha$ completes the proof.$\Box$


\vskip 1cm

{\noindent
\parbox[t]{65truemm}
{Hans Havlicek\\
Abteilung f\"ur Lineare Algebra\\
und Geometrie\\
Technische Universit\"at\\
Wiedner Hauptstra\ss e 8--10\\
A-1040 Wien, Austria}
\hfill
\parbox[t]{65truemm}
{Corrado Zanella\\
Dip.\ di Matematica Pura\\
ed Applicata\\
Universit\`a di Padova\\
via Belzoni 7\\
I-35131 Padova, Italy}}

\end{document}